\documentclass[journal]{IEEEtran}

\usepackage{cite}

\usepackage[cmex10]{amsmath}
\usepackage[pdftex]{graphicx}
\usepackage{amsfonts}
\usepackage{mathrsfs}
\usepackage{relsize}
\usepackage{mathtools} 
\usepackage{subcaption}
\DeclareMathOperator{\tr}{tr} 

\begin{document}
%
% paper title
% can use linebreaks \\ within to get better formatting as desired
\title{Model-based Optimization of Compressive Antennas for High-Sensing-Capacity Applications}

\author{\IEEEauthorblockN{Richard Obermeier and Jose~Angel~Martinez-Lorenzo}\\% <-this % stops a space
\IEEEauthorblockA{%line 1: dept. name (if applicable)\\
Northeastern University \\%line 2: name of organization, acronyms acceptable\\
Boston, MA, USA\\%line 3: City, State/Province, Country\\
obermeier.ri@husky.neu.edu, jmartinez@coe.neu.edu}}%line 4: e-mail address if desired\\}}

\maketitle

\noindent \begin{abstract}
This paper presents a novel, model-based compressive antenna design method for high sensing capacity imaging applications. Given a set of design constraints, the method maximizes the sensing capacity of the compressive antenna by varying the constitutive properties of scatterers distributed along the antenna. Preliminary 2D design results demonstrate the new method's ability to produce antenna configurations with enhanced imaging capabilities. 
%Given a set of parameter constraints, such as transmitting antenna positions, imaging positions, and frequencies, the design method attempts to maximize the channel capacity of the dyadic Green's function matrix. Numerical results show that the proposed design method significantly improves the channel capacity of the Green's function matrix,  thereby improving the antenna's imaging capabilities. 
\end{abstract}


%{\let\newpage\relax\maketitle}
\section{Introduction}
Sensing systems attempt to extract as much information as possible about an object or region of interest by recording a set of independent measurements. The number of measurements, and the degree of their independence, determine how much information a sensing system can extract. Recent papers \cite{Martinez2015a,Heredia2015} have introduced the concept of a compressive reflector antenna for use in millimeter wave imaging applications. The compressive reflector antenna operates in a manner similar to that of the coded apertures utilized in optical imaging applications \cite{busboom1997coded,de2009sub,marcia2008compressive}: by introducing scatterers to the surface of a traditional reflector antenna, the compressive antenna encodes a pseudo-random phase front on the scattered electric field. By modifying the encoded wavefront from measurement to measurement, for example by rotating the reflector or by electrically changing the constitutive properties of the scatterers, compressive sensing techniques \cite{Candes2006,Donoho2006,massa2015} can be employed with improved performance over the traditional reflector antenna. This paper describes a numerical method for optimizing the constitutive parameters of the scattering elements in order to design compressive reflector antennas with high sensing capacity. %Our design approach uses a numerical model based upon finite differences in the frequency domain (FDFD) \cite{Rappaport2001} in order to maximize the channel capacity of the dyadic Green's function matrix. 

The remainder of this paper is organized as follows. Section \ref{sec:motivation} shows why the sensing, or channel capacity is a reasonable metric to use when assessing the sensing capabilities of an antenna. Section \ref{sec:design} describes a general compressive reflector antenna design approach, which optimizes the constitutive parameters of the scattering elements under realistic constraints in order to maximize the sensing capacity. Section \ref{sec:simple} describes a simplified design approach, which optimizes over the constitutive parameters of the reflector elements under box constraints. Section \ref{sec:results} presents preliminary 2D antenna design results, which demonstrate the new method's ability to increase the channel capacity. Finally, the conclusions are presented in Section \ref{sec:conclusions}.
%Our design approach utilizes a numerical model based upon finite differences in the frequency domain (FDFD) \cite{Rappaport2001}  

\section{Motivation}
\label{sec:motivation}
One of the key features of next generation sensing and imaging systems will be the ability to maximize the sensing capacity \cite{Martinez2015a} that is, the information transfer efficiency between the pixels in the imaging region and the data measured by the system. This can occur when the mutual information of successive measurements is as low as possible. One way to achieve this goal is to dynamically control the ``wavefield information'' that is encoded in several dimensions (i.e. time, space, frequency, phase, polarization and novel angular momentum) in any given sensing experiment $-$ a methodology known as multi-dimensional codification. 

The sensing matrix $\mathbf{A}$ of a compressive antenna working in a monostatic configuration is fully characterized by its radiation pattern. This radiation pattern can be derived from the dyadic Green's functions $\mathbf{G}(\mathbf{r}, \mathbf{r}',\omega)$ of the compressive antenna, which can be expressed as the solution to the vector wave equation:
\begin{align}
\nabla\times \frac{1}{\mu(\mathbf{r},\omega)}\nabla \times \mathbf{G}(\mathbf{r},\mathbf{r}',\omega) &- \omega^2\epsilon(\mathbf{r},\omega)\mathbf{G}(\mathbf{r}, \mathbf{r}', \omega) \nonumber \\
&= \mathbf{\tilde{I}} \delta\left(\mathbf{r}-\mathbf{r}'\right) \label{eq:green_prob}
\end{align}
Using this definition for $\mathbf{G}(\mathbf{r}, \mathbf{r}',\omega)$, the electric field radiated by a compressive antenna is given by:
\begin{align}
\mathbf{E}(\mathbf{r}, \omega) &= \jmath\omega \int \mathbf{G}(\mathbf{r}, \mathbf{r}', \omega) \cdot \mathbf{I}(\mathbf{r}', \omega) d\mathbf{r}' \label{eq:green_sol}
\end{align}
where $\mathbf{I}(\mathbf{r}', \omega)$ describe the sources used to excite the reflector. Eq. \ref{eq:green_sol} can be discretized into a linear system of equations, $\mathbf{E}_\omega = \mathbf{G}_\omega\mathbf{S}_\omega$, where the constants have been absorbed into the source term $\mathbf{S}_\omega$. Due to the reciprocity theorem, this relationship can also be used to described the electric fields scattered from an object of interest. In this case, the source distribution $\mathbf{S}_\omega$ can be interpreted as the set of ``contrast sources'' located within the imaging region \cite{VanDenBerg2001}, and the field vector $\mathbf{E}_\omega$ are the scattered electric fields measured by the receivers. The sensing capabilities of the compressive antenna can thus be obtained by analyzing the Green's function matrix through its singular value decomposition, $\mathbf{G}_\omega = \mathbf{U}\mathbf{\Sigma}\mathbf{V}^H$, where $\mathbf{U}$ and $\mathbf{V}$ are orthonormal matrices and $\mathbf{\Sigma}$ is a rectangular diagonal matrix with diagonal elements $\sigma_i \ge 0$; the $\sigma_i$ are known as the singular values. When a ``contrast source'' distribution $\mathbf{v}_i$ excites the system, an electric field $\sigma_i\mathbf{u}_i$ is generated at the receivers. If the singular values are poorly conditioned, i.e. $\sigma_{\text{max}}/\sigma_{\text{min}} \gg 1$, then it becomes difficult to distinguish different source distributions from each other using only the measured fields. The sensing, or channel capacity of the Green's function matrix can be used to compactly assess the imaging capabilities of the compressive antenna. For high Signal to Noise Ratios (SNR), the sensing capacity can be expressed as \cite{Proakis2008}:
%\begin{align}
%C &= \sum_{t=1}^{T}\log_2\left(1 + \frac{P_t}{TN_0}\sigma_t^2\right) \nonumber \\
%&\approx T\log_2\left(\sum_{t=1}^T \frac{P_t}{N_0}\right) + \sum_{t=1}^T \log_2\left(\frac{\sigma_t^2}{T}\right) \label{eq:capacity}
%\end{align}
\begin{align}
C & \approx T\log_2\left(\sum_{t=1}^T \frac{P_t}{N_0}\right) + \sum_{t=1}^T \log_2\left(\frac{\sigma_t^2}{T}\right) \label{eq:capacity}
\end{align}
The proposed antenna design method seeks to maximize the sensing capacity of the Green's function matrix by optimizing the constitutive parameters $\mu(\mathbf{r},\omega)$ and $\epsilon(\mathbf{r},\omega)$ of the scattering elements placed upon the  reflector. Since the channel capacity penalizes over the logarithm of the singular values, this approach favors systems with small condition numbers $\sigma_{\text{max}}/\sigma_{\text{min}}$. 

\section{A General Design Approach}
\label{sec:design}%  $\{\mathbf{r}_1^t, \mathbf{r}_2^t, \ldots, \mathbf{r}_T^t\}$  $\{\mathbf{r}_1^i, \mathbf{r}_2^i, \ldots, \mathbf{r}_M^i\}$ $\{\omega_1, \omega_2, \ldots, \omega_K\}$ $\{\mathbf{r}_1^s, \mathbf{r}_2^s, \ldots, \mathbf{r}_N^s\}$ 
In the optimization problem, the transmitting antenna system is described by a set of current sources located at $T$ locations. Each transmitting antenna excites the $M$ positions in the imaging region with stepped-frequency waveforms at $K$ frequencies. The design procedure optimizes the constitutive properties $\epsilon(\mathbf{r},\omega)$ and $\mu(\mathbf{r},\omega)$ of scattering elements located at $N$ positions along the reflector. In order to allow the scattering elements to be dispersive, the permittivity and permeability of the scatterers at the $k-$th frequency will be jointly represented by the variable $\mathbf{x}_k$. With this convention, the matrix $\mathbf{G}_{k}\left(\mathbf{x}_k\right) \in \mathbb{C}^{3 M \times 3 T}$ can be defined as the Green's function matrix for sources radiating at frequency $\omega_k$, located at the $T$ transmitter positions, and evaluated at the $M$ positions in the imaging region. This matrix is a nonlinear function of the design variables $\mathbf{x}_k$. By concatenating the Green's function matrices for multiple frequencies, the multi-frequency Green's function matrix $\mathbf{G}(\mathbf{x}) \in \mathbb{C}^{3 M \times 3 K T}$ can be expressed as:
\begin{align}
\mathbf{G}(\mathbf{x}) &= \mathbf{G}(\mathbf{x}_1, \mathbf{x}_2, \ldots, \mathbf{x}_K) \nonumber \\
&= \big[\mathbf{G}_{1}\left(\mathbf{x}_1\right), \mathbf{G}_{2}\left(\mathbf{x}_2\right), \ldots, \mathbf{G}_{K}\left(\mathbf{x}_K\right)\big] \label{eq:bigG}
\end{align}
where the vector $\mathbf{x}$ is the vector of concatenated design variables for each frequency. Assuming that $M > KT$, the channel capacity maximization problem can be expressed as a non-convex ``max-det'' problem:
 \begin{align}
&\text{maximize} ~~\log\det\left(\mathbf{G}^H(\mathbf{x}) \mathbf{G}(\mathbf{x})\right) \label{eq:maxdet} \\
&\text{subject to} ~~h_q(\mathbf{x}) \le 0, ~~~~ q = 1, \cdots, Q \nonumber  \nonumber \\
&~~~~~~~~~~~~~~c_p(\mathbf{x}) = 0, ~~~~ p = 1, \cdots, P \nonumber
%&1 &\mathbf{z} \ge \mathbf{0}
\end{align}
It is easy to show that maximizing the log-determinant of the Grammian matrix is equivalent to maximizing the channel capacity. Since $\mathbf{G}^H(\mathbf{x}) \mathbf{G}(\mathbf{x}) = \mathbf{V}(\mathbf{x})\mathbf{\Sigma}^2(\mathbf{x})\mathbf{V}^H(\mathbf{x})$, $\log\det\left(\mathbf{G}^H(\mathbf{x}) \mathbf{G}(\mathbf{x})\right) = \sum_{t=1}^T \log\left(\sigma_t^2\right)$, which differs from the channel capacity defined in Eq. \ref{eq:capacity} only by constants.
%Since $\mathbf{G}^H(\mathbf{x}) \mathbf{G}(\mathbf{x}) = \mathbf{V}(\mathbf{x})\mathbf{\Sigma}^2(\mathbf{x})\mathbf{V}^H(\mathbf{x})$, $\log\det\left(\mathbf{G}^H(\mathbf{x}) \mathbf{G}(\mathbf{x})\right) = \sum_{t=1}^T \log\left(\sigma_t^2\right)$, and so this optimization problem maximizes the channel capacity.

The constraint functions $h_q(\mathbf{x})$ and $c_p(\mathbf{x})$ can be non-convex and depend upon the specific design constraints placed on the dielectric scatterers. For example, if the scatterers are restricted to non-dispersive materials, then the equality constraint functions force the design variables $\mathbf{x}_1, \mathbf{x}_2, \ldots, \mathbf{x}_K$ to produce the same permittivity and conductivity. As another example, if metamaterial scattering elements are disallowed, then the inequality constraint functions force the design variables to produce dielectric constants $\ge 1$. 

\section{A Simplified Design Approach}
\label{sec:simple}
This section describes how to solve a simplified version of Eq. \ref{eq:maxdet}. In this approach, both the scatterers and the background medium at the scatterer locations are assumed to be non-dispersive and non-conductive, so that the design variables $\mathbf{x}_1, \mathbf{x}_2, \ldots, \mathbf{x}_K$ are equal and are real-valued. Moreover, the constraints simply restrict the electric permittivities and magnetic permeabilities of the scatterers to lie within specified ranges, $[\epsilon_L, \epsilon_R]$ and $[\mu_L, \mu_R]$. The simplified optimization problem can therefore be expressed as:
 \begin{align}
&\text{maximize} ~~\log\det\left(\mathbf{G}^H(\mathbf{x}) \mathbf{G}(\mathbf{x})\right) \label{eq:maxdet2} \\
&\text{subject to} ~~\mathbf{x}_L \le \mathbf{x} \le \mathbf{x}_R \nonumber
%&1 &\mathbf{z} \ge \mathbf{0}
\end{align}
Eq. \ref{eq:maxdet2} can be solved efficiently using the nonlinear conjugate gradient method \cite{nocedal2006numerical}. This method requires expressions for the gradient of the cost function $\log\det\mathbf{F}(\mathbf{x}) = \log\det\left(\mathbf{G}^H(\mathbf{x}) \mathbf{G}(\mathbf{x})\right)$. Assuming that $\mathbf{F}(\mathbf{x})$ that is invertible, the partial derivatives $\frac{\partial}{\partial x_l}\log\det\mathbf{F}(\mathbf{x})$ and $\frac{\partial \mathbf{F}(\mathbf{x}) }{\partial x_l}$ are:
\begin{align}
\frac{\partial}{\partial x_l} \log \det \mathbf{F}(\mathbf{x}) &= \tr\left(\mathbf{F}^{-1}(\mathbf{x})\frac{\partial \mathbf{F}(\mathbf{x})}{\partial x_l} \right) \\
\frac{\partial \mathbf{F}(\mathbf{x}) }{\partial x_l} &=  \left(\frac{\partial \mathbf{G}(\mathbf{x})}{\partial x_l}\right)^H\mathbf{G}(\mathbf{x}) + \mathbf{G}^H(\mathbf{x})\frac{\partial \mathbf{G}(\mathbf{x})}{\partial x_l}
\end{align} 
A close examination of Eq. \ref{eq:bigG} reveals hat the partial derivatives $\frac{\partial\mathbf{G}(\mathbf{x})}{\partial x_l}$ consist of the partial derivatives $\frac{\partial\mathbf{G}_k(\mathbf{x})}{\partial x_l}$. By defining $\mathbf{H}_k(\mathbf{x})$ as the discretized version of the Helmholtz operator for frequency $k$, the Green's function matrix $\mathbf{G}_k(\mathbf{x})$ can be expressed as:
\begin{align}
\mathbf{G}_k(\mathbf{x}) = \mathbf{\Phi} \mathbf{H}_k^{-1}(\mathbf{x}) \mathbf{\Psi}
\end{align}
where $\mathbf{\Phi} \in \mathbb{C}^{3 M \times 3 L}$, $\mathbf{H}_k(\mathbf{x}) \in \mathbb{C}^{3 L \times 3 L}$, and $\mathbf{\Psi} \in \mathbb{C}^{3 L \times 3 T}$. The matrices $\Phi$ and $\Psi$ are subsampling matrices corresponding to the imaging and transmitter positions respectively. From this relationship, the partial derivatives $\frac{\partial\mathbf{G}_k(\mathbf{x})}{\partial x_l}$ take the following form:
\begin{align}
\frac{\partial\mathbf{G}_k(\mathbf{x})}{\partial x_l} = -\mathbf{\Phi}\mathbf{H}_k^{-1}(\mathbf{x})\frac{\partial \mathbf{H}_k(\mathbf{x})}{\partial x_l} \mathbf{H}_k^{-1}(\mathbf{x})\mathbf{\Psi} \label{eq:partialGreen}
\end{align}
The elements of the partial derivative matrix $\frac{\partial \mathbf{H}_k(\mathbf{x})}{\partial x_l}$ differ depending upon whether $x_l$ is permittivity or permeability. If $x_l$ is the permittivity $\epsilon_j$ at position $j$, then the partial derivative matrix takes the form:
%\begin{equation}
%\left(\frac{\partial \mathbf{H}_k(\mathbf{x})}{\partial \epsilon_p}\right)_{mn} = 
%\begin{cases}
%\omega_k^2& m = n =p \\%\left(k_{j}^{(b,l)}\right)^2 & m = n = l \\
%0 & \text{otherwise}
%\end{cases}
%\end{equation}
\begin{equation}
\frac{\partial \mathbf{H}_k(\mathbf{x})}{\partial \epsilon_j} = \omega_k^2 \operatorname{diag}(\mathbf{1}_3 \otimes \boldsymbol{\delta}_{ij})
\end{equation}
where $\otimes$ is the Kronecker product and $\boldsymbol{\delta}_{ij} \in \mathbb{C}^{L}$ is the Kronecker delta function expressed as a vector, i.e. the $j-th$ element of $\boldsymbol{\delta}_{ij}$ equals one and all others equal zero. % the elements of the vector $\tilde{\boldsymbol{\epsilon}}_p$ are:
%\begin{equation}
%\left(\tilde{\boldsymbol{\epsilon}}_p\right)_m = 
%\begin{cases}
%\omega_k^2& m =p \\%\left(k_{j}^{(b,l)}\right)^2 & m = n = l \\
%0 & \text{otherwise}
%\end{cases}
%\end{equation}
%where the elements of the partial derivative matrix $\frac{\partial \mathbf{H}_j(\mathbf{x})}{\partial x_l}$ are given by:
%\begin{equation}
%\left(\frac{\partial \mathbf{H}_j(\mathbf{x})}{\partial x_l}\right)_{mn} = 
%\begin{cases}
%\omega_j^2(\boldsymbol{\mu}_b)_l (\boldsymbol{\epsilon}_b)_l & m = n = l \\%\left(k_{j}^{(b,l)}\right)^2 & m = n = l \\
%0 & \text{otherwise}
%\end{cases}
%\end{equation}
If $x_l$ is the permeability $\mu_j$ at position $j$, then the partial derivative matrix takes the form:
\begin{equation}
\frac{\partial \mathbf{H}_k(\mathbf{x})}{\partial \mu_j}= -\frac{1}{\mu_j^2}\mathbf{L}_c \operatorname{diag}(\mathbf{1}_3 \otimes \boldsymbol{\delta}_{ij})\mathbf{L}_c
\end{equation}
where $\mathbf{L}_c$ is the discretized curl operator.
% and the elements of the vector $\tilde{\boldsymbol{\mu}}_p$ are:
%\begin{equation}
%\left(\tilde{\boldsymbol{\mu}}_p\right)_m = 
%\begin{cases}
%\frac{1}{\mu_p}& m =p \\%\left(k_{j}^{(b,l)}\right)^2 & m = n = l \\
%0 & \text{otherwise}
%\end{cases}
%\end{equation}
Computation of these derivatives requires $K(N+T)$ calls to a forward model solver at each iteration in order to compute the necessary Green's functions. %In our work, we use a forward solver based on finite differences in the frequency domain \cite{Rappaport2001}. 

\section{Results}
\label{sec:results}
This sections presents preliminary antenna design results, which were generated using the simplified algorithm and a 2D forward model solver based on finite differences in the frequency domain (FDFD) \cite{Rappaport2001}. The design method was executed for two different antenna configurations, where the antenna operated in reflection mode and transmission mode. In reflection mode, dielectric scatterers are added to the surface of a Perfect Electric Conductor (PEC) reflector in order to further perturb the fields scattered by the reflector. In transmission mode, dielectric scatterers are placed in between the transmitting antennas and the imaging region in order to perturb the fields that otherwise radiate in the homogeneous background medium.

\begin{figure}[h]
\begin{subfigure}[h]{1\columnwidth}
\centering
\includegraphics[width=5.5cm, clip=true]{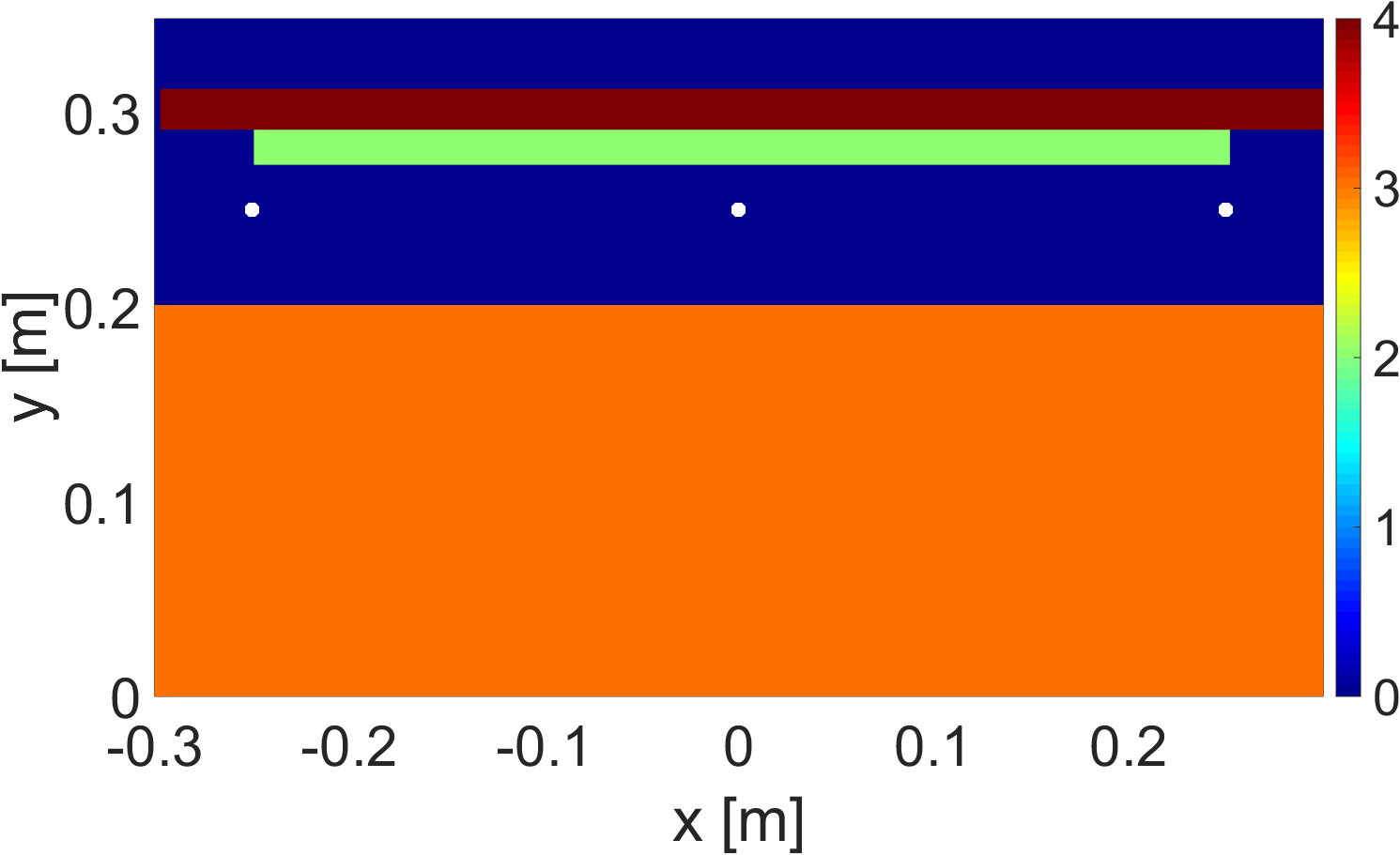}
\caption{}
\label{fig:rx_config}
\end{subfigure}
%\hskip 0.05\columnwidth
\begin{subfigure}[h]{1\columnwidth}
\centering
\includegraphics[width=5.5cm, clip=true]{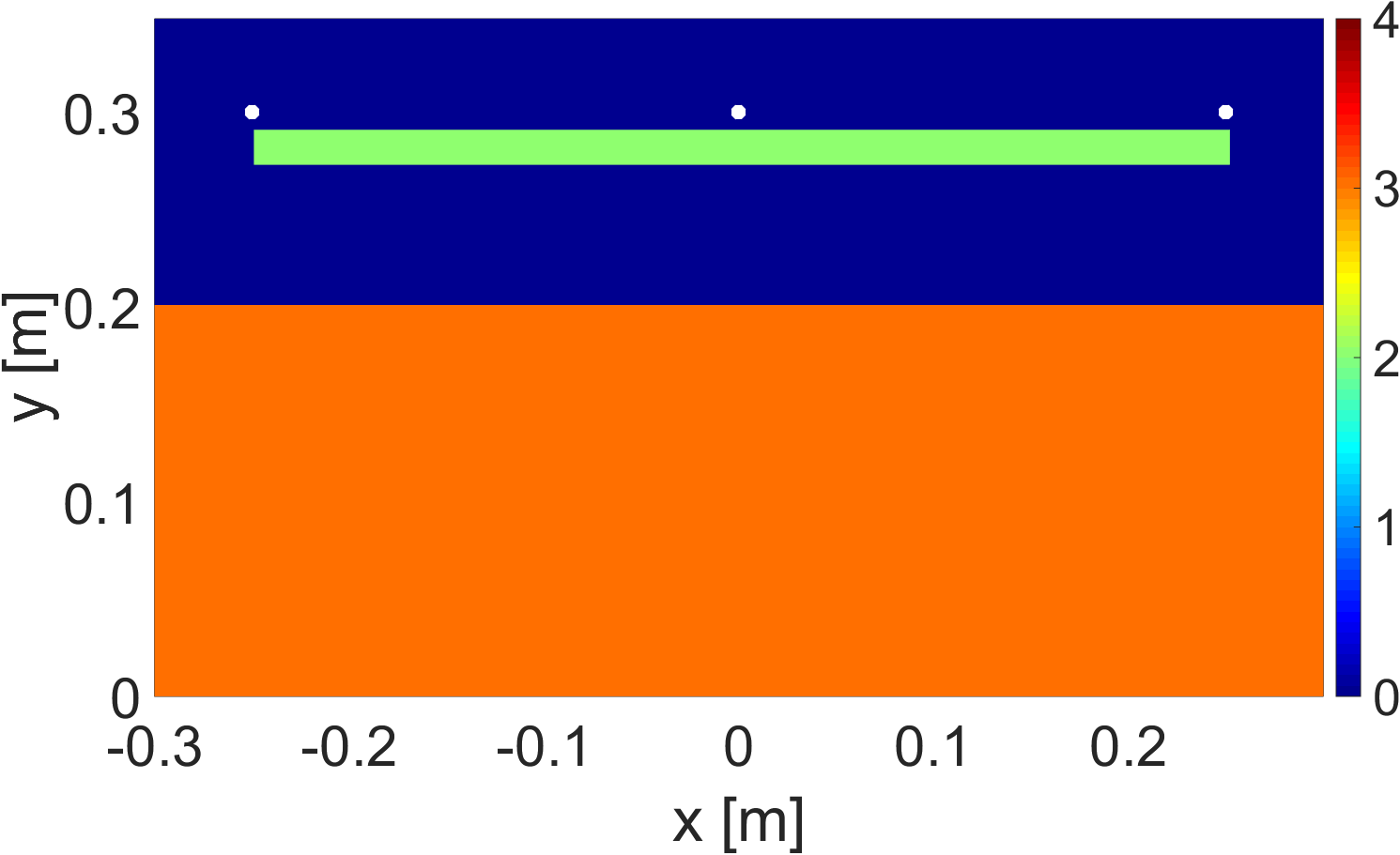}
\caption{}
%\caption{\small Configuration for the transmission mode design problem. Light Blue = Transmitter locations, Orange - Imaging region, Green = Scatterer locations.}
\label{fig:tx_config}
\end{subfigure}
\caption{\small Configuration for the compressive antenna operating in (a) reflector mode and (b) transmission mode. Light Blue = Transmitter locations, Orange - Imaging region, Green = Scatterer locations, Red = PEC.}
\end{figure}

Figures \ref{fig:rx_config} and \ref{fig:tx_config} display the design configurations for the reflection and transmission mode problems respectively. In both modes, three line source antennas, represented by the white circles, were used to excite the free-space imaging region, colored in orange. The green pixels represent the locations of the scatterers to be optimized, and the red pixels in the reflection mode configuration represent the PEC. The antennas were constrained to transmit at five frequencies linearly spaced between $3.1$GHz and $3.5$GHz, and the dielectric constant of the scatterers was constrained to the range $[1,10]$; the magnetic permeability was restricted to $\mu = \mu_0$. 

%\begin{figure}[h!]
%\begin{minipage}{.47\linewidth}
%\centering
%\includegraphics[height=3.2cm, clip=true]{./reflector_design_config_mod.png}
%\caption{\small Configuration for the reflector mode design problem. Light Blue = Transmitter locations, Orange - Imaging region, Green = Scatterer locations, Red = PEC.}
%\label{fig:rx_config}
%\end{minipage}
%\hskip .03\linewidth
%\begin{minipage}{.47\linewidth}
%\centering
%\includegraphics[height=3.2cm, clip=true]{./tx_design_config_mod.png}
%\caption{\small Configuration for the transmission mode design problem. Light Blue = Transmitter locations, Orange - Imaging region, Green = Scatterer locations.}
%\label{fig:tx_config}
%\end{minipage}
%\end{figure}

\begin{figure}[h]
\begin{subfigure}[h]{1\columnwidth}
\centering
\includegraphics[width=5.9cm, clip=true]{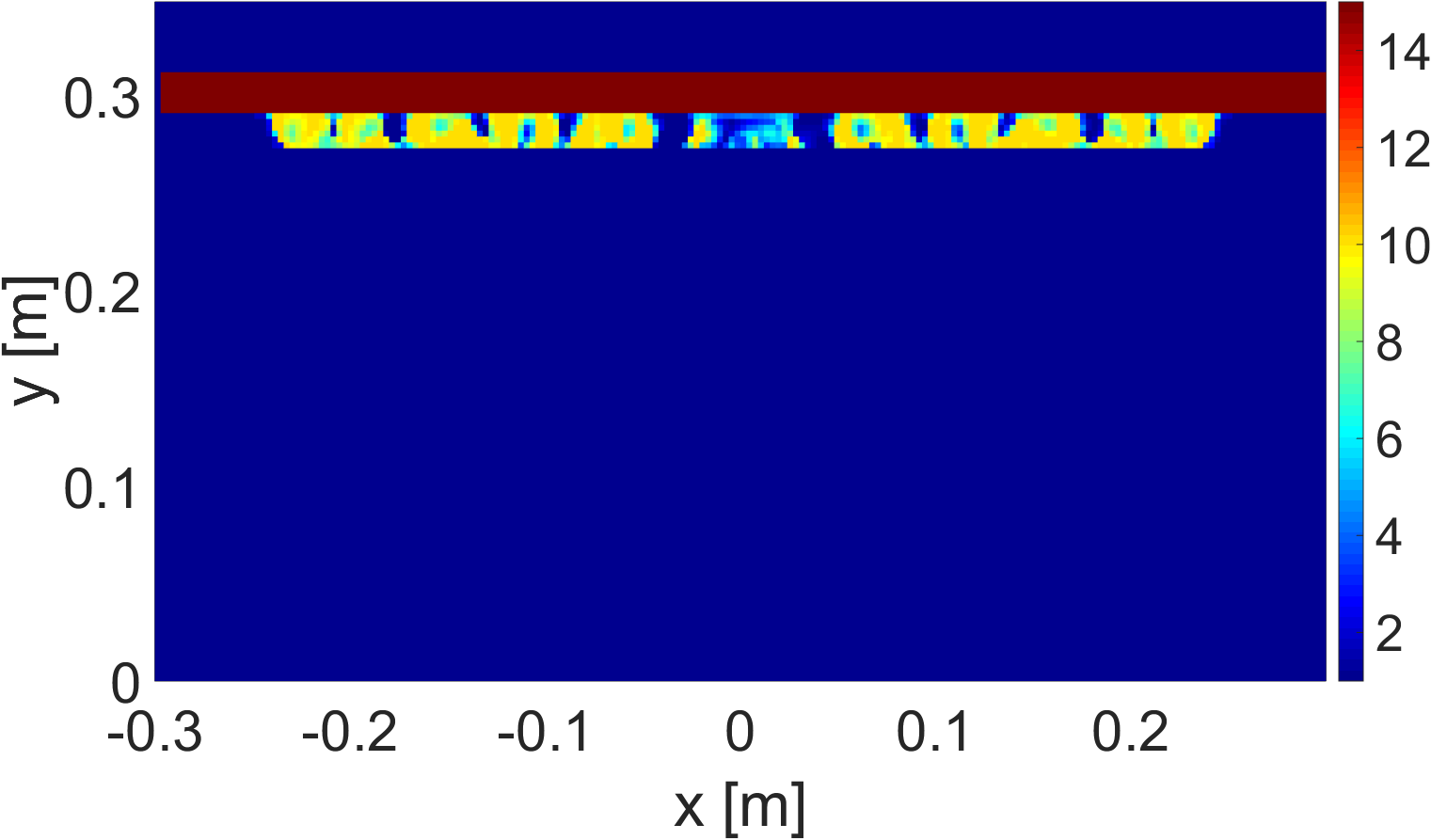}
\caption{}
\label{fig:rx_permittivity}
\end{subfigure}
%\hskip 0.05\columnwidth
\begin{subfigure}[h]{1\columnwidth}
\centering
\includegraphics[width=5.9cm, clip=true]{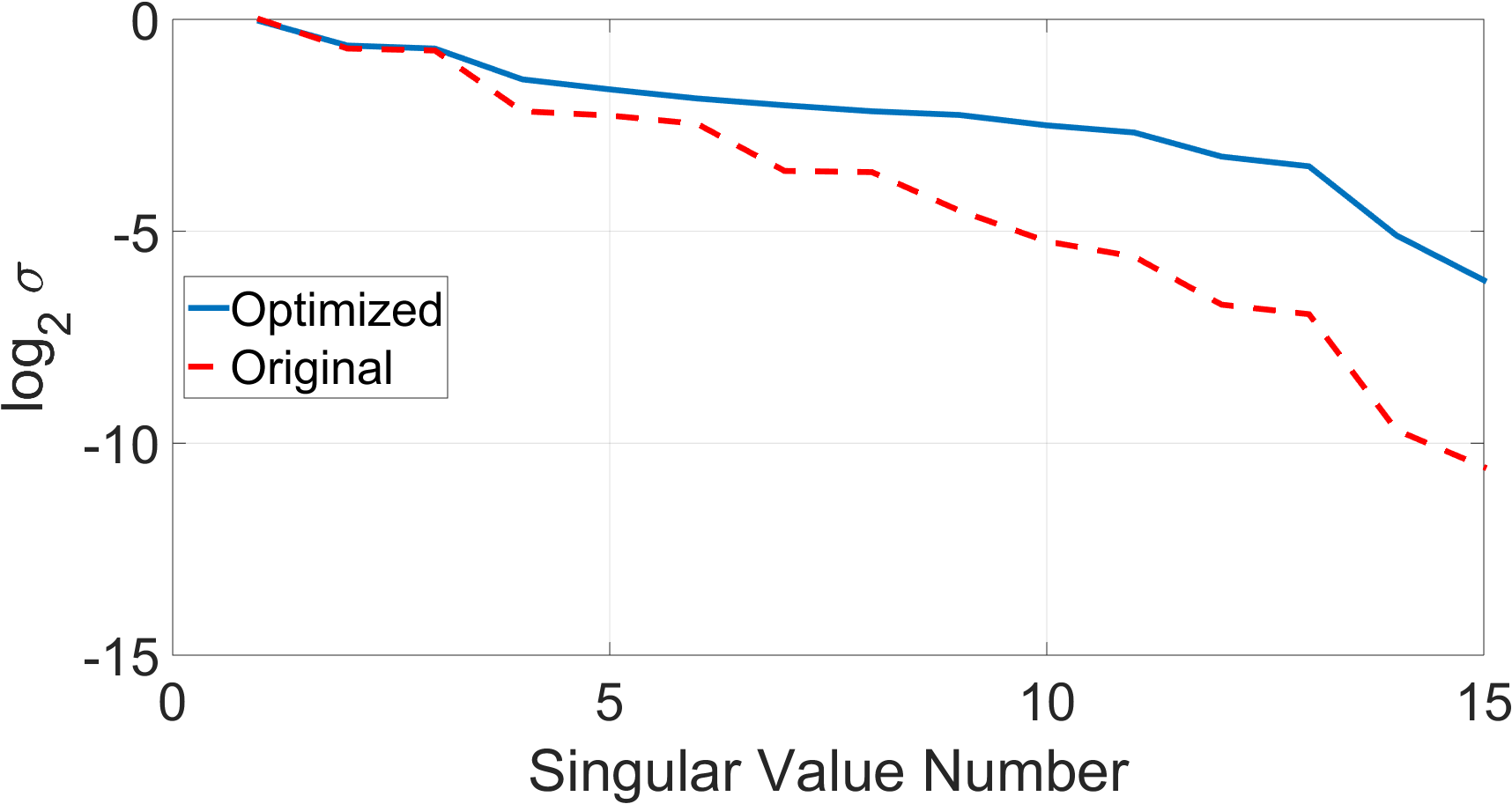}
\caption{}
\label{fig:rx_log2}
\end{subfigure}
\caption{\small Reflection mode configuration: (a) optimized permittivity. (b) $\log_2$ of the singular values of Green's function matrices.}
\end{figure}

\begin{figure}[b!]
\begin{subfigure}[h]{\columnwidth}
\centering
\includegraphics[width=5.9cm, clip=true]{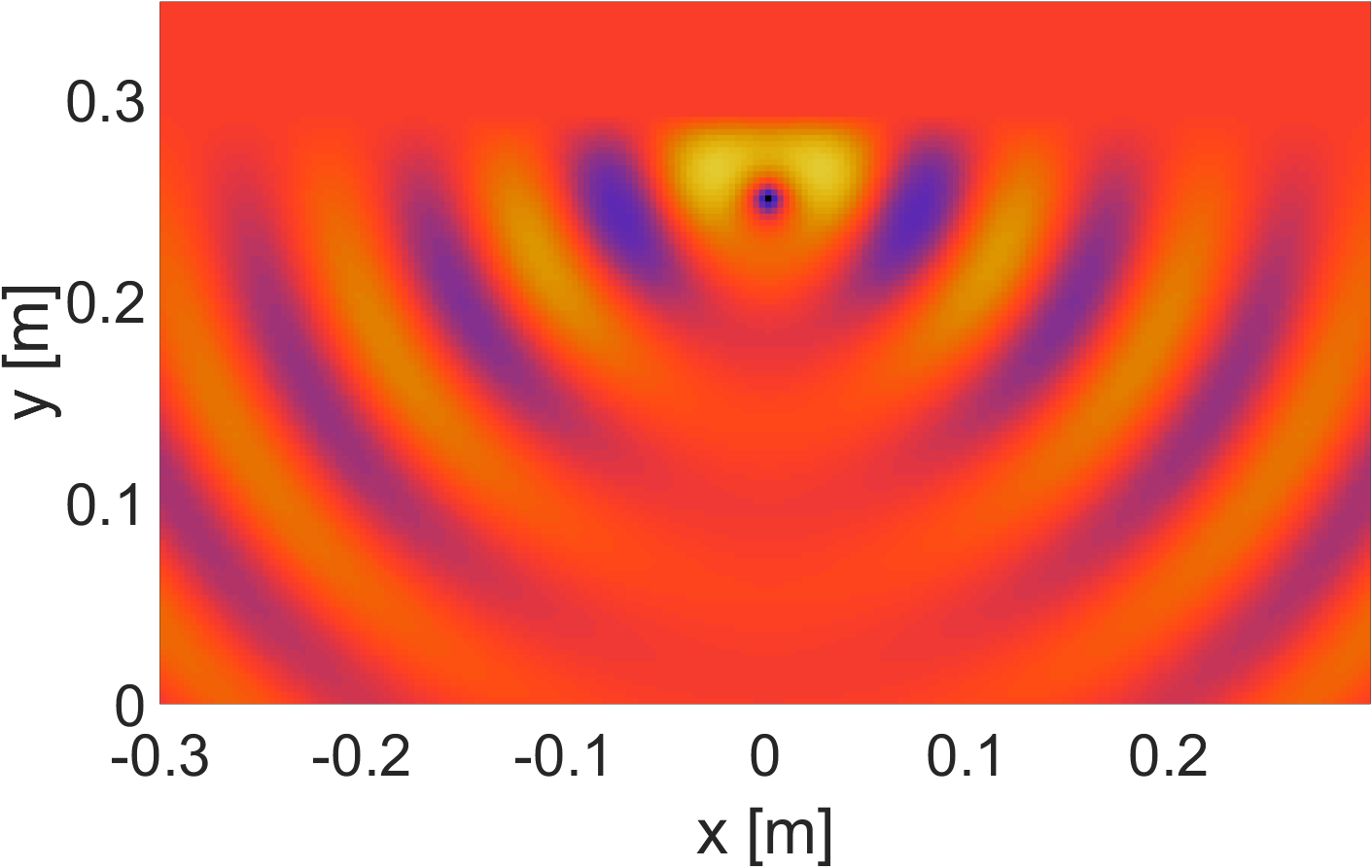}
\caption{}
\label{fig:rx_green_orig}
\end{subfigure}
\begin{subfigure}[h]{\columnwidth}
\centering
\includegraphics[width=5.9cm, clip=true]{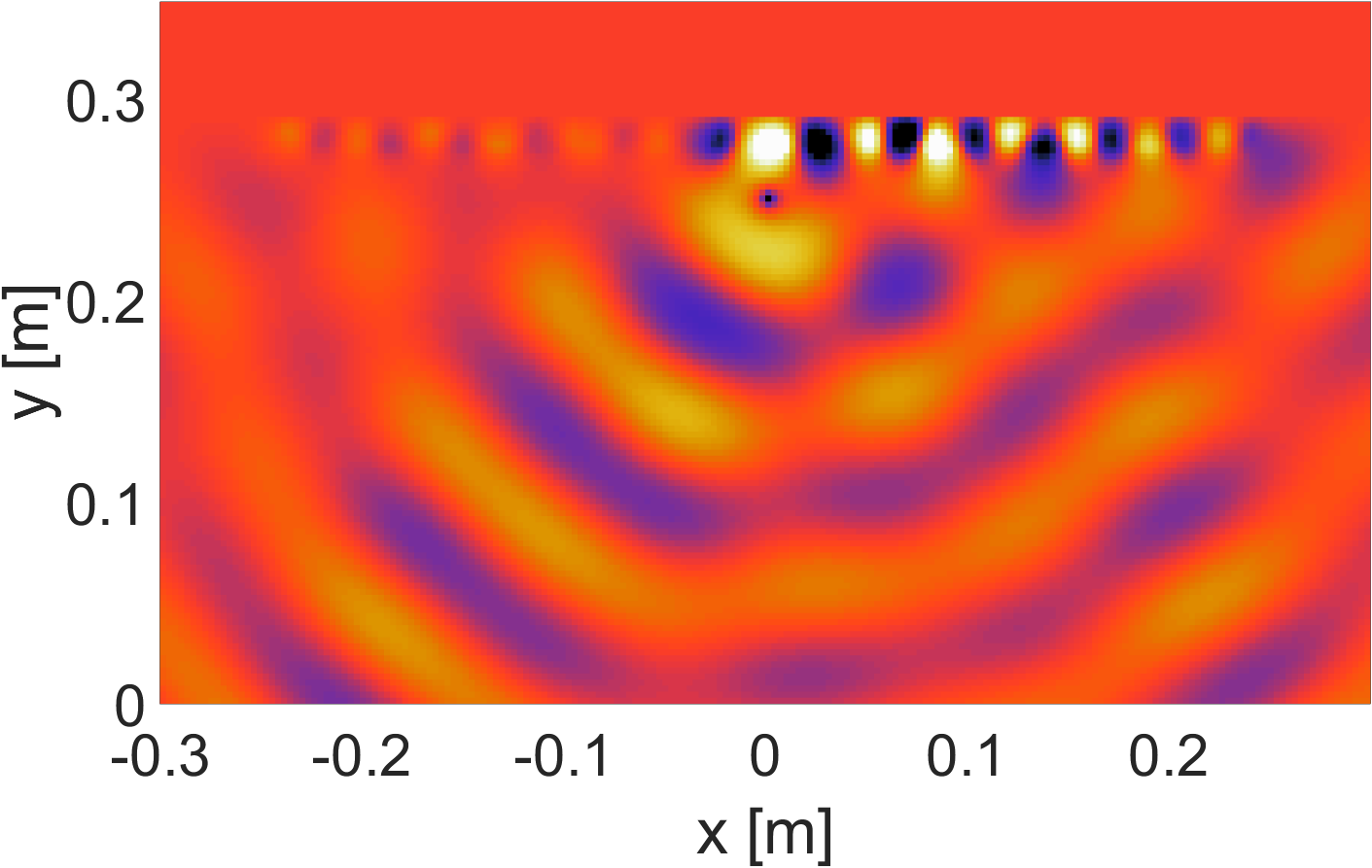}
\caption{}
\label{fig:rx_green_opt}
\end{subfigure}
\caption{\small Green's functions of the compressive antenna operating in reflection mode at $3.5$GHz with (a) the original design and (b) optimized design.}
\end{figure}

Figure \ref{fig:rx_permittivity} displays the optimized permittivity distribution for the reflection mode problem. It is important to note that the design problem of Eq. \ref{eq:maxdet2} is non-convex, and so it is probable that the solution displayed in Figure \ref{fig:rx_permittivity} is only a locally optimal solution. In practice, the optimization problem can be solved several times using different starting points until an antenna design with suitable sensing capacity is found. Figure \ref{fig:rx_log2} displays the $\log_2$ of the singular values of the Green's function matrices for the original and optimized antennas. In this configuration, the optimized design increases the channel capacity by $29$ bits per pixel and decreases the condition number by a factor of $22$, from approximately $1510$ to $69$. Figures \ref{fig:rx_green_orig} and \ref{fig:rx_green_opt} display the Green's function of the middle antenna radiating at $3.5$GHz for the original and optimized antennas respectively. While the phase front of the original antenna is fairly uniform as a function of distance from the transmitter, the phase front of the optimized antenna is noticeably perturbed.   

\begin{figure}[b!]
\begin{subfigure}[h]{1\columnwidth}
\centering
\includegraphics[width=5.9cm, clip=true]{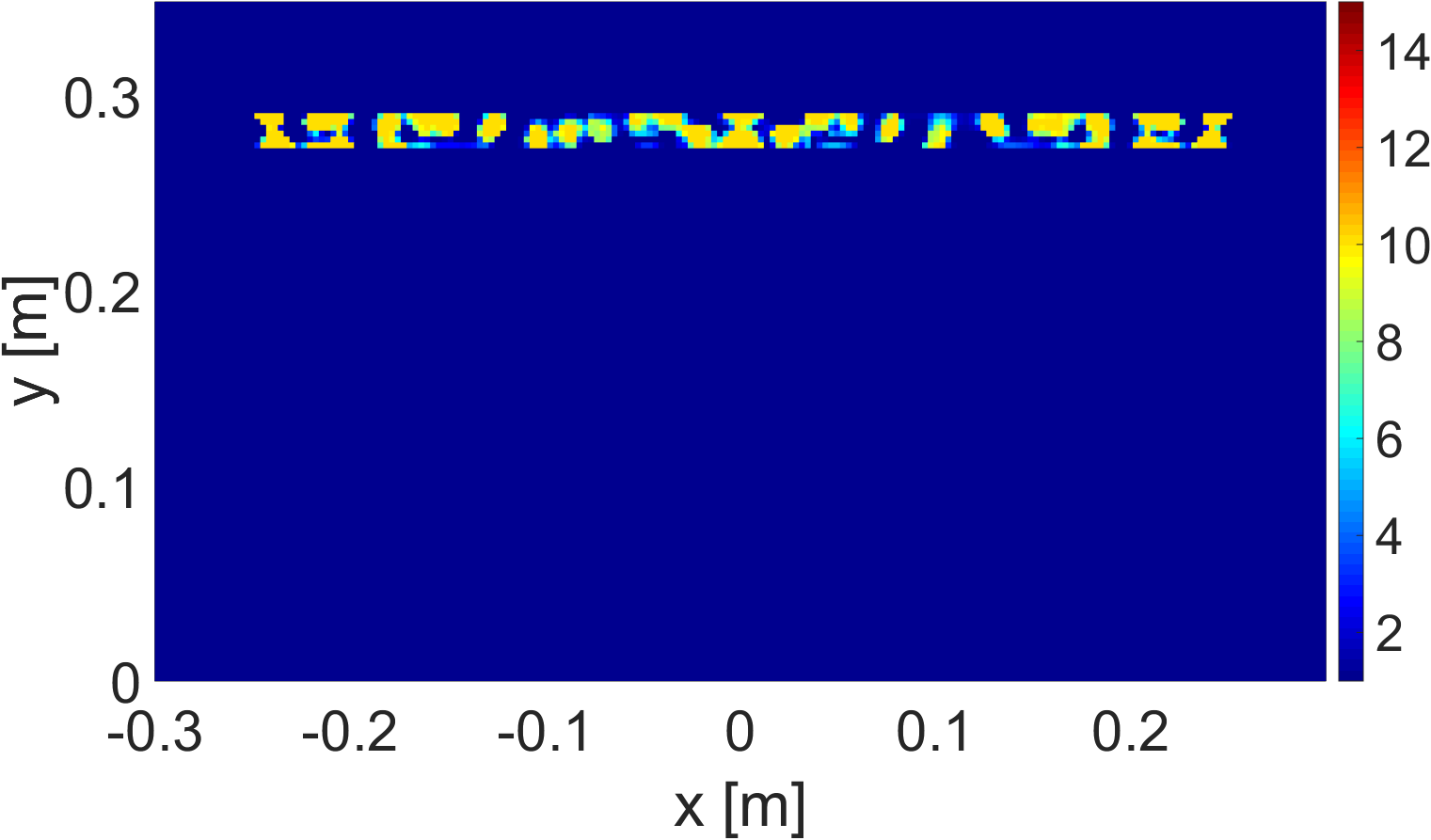}
\caption{}
\label{fig:tx_permittivity}
\end{subfigure}
%\hskip 0.05\columnwidth
\begin{subfigure}[h]{1\columnwidth}
\centering
\includegraphics[width=5.9cm, clip=true]{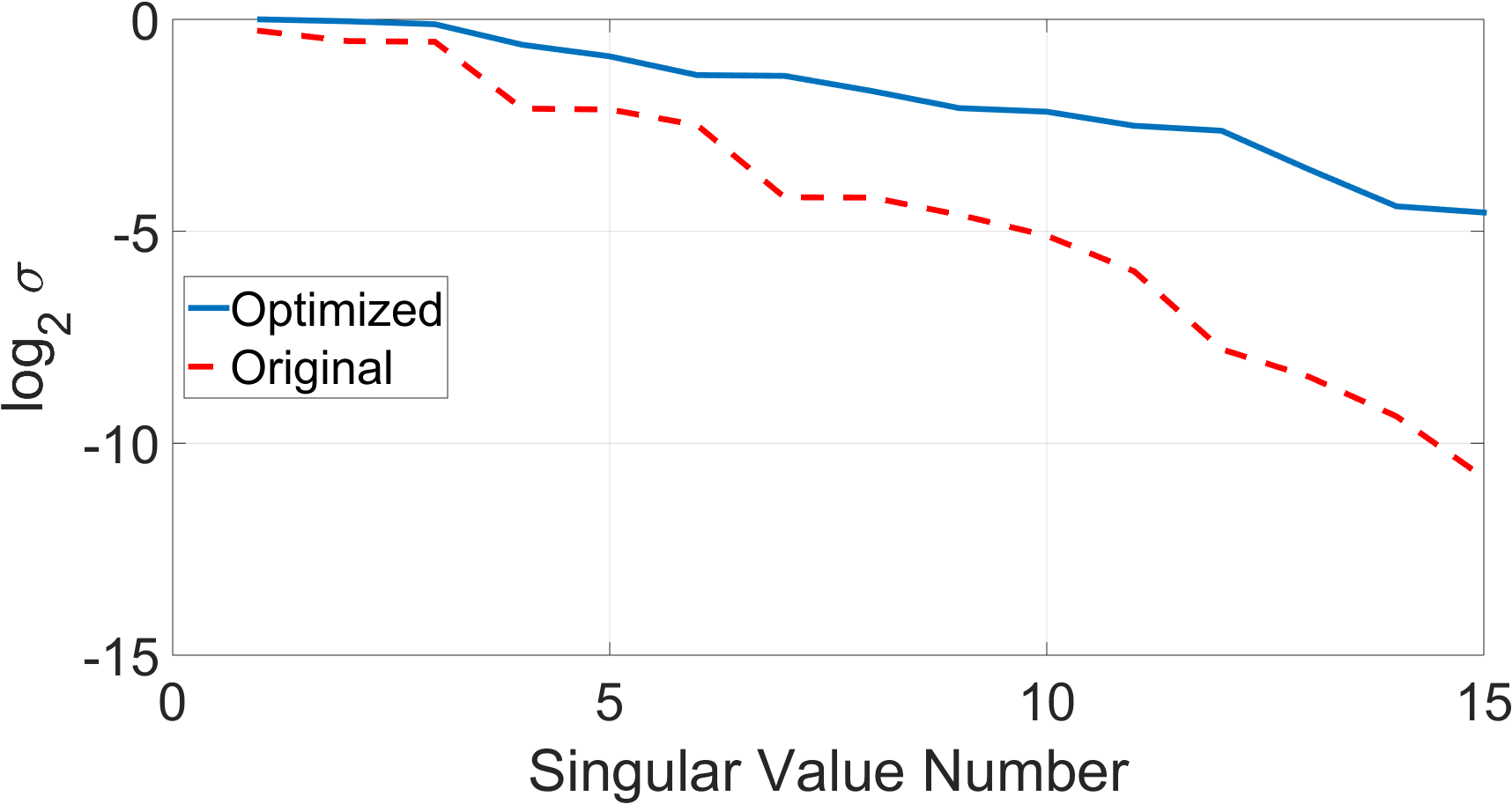}
\caption{}
\label{fig:tx_log2}
\end{subfigure}
\caption{\small Transmission mode configuration: (a) optimized permittivity. (b) $\log_2$ of the singular values of Green's function matrices.}
\end{figure}

\begin{figure}[t]
\begin{subfigure}[h]{\columnwidth}
\centering
\includegraphics[width=5.9cm, clip=true]{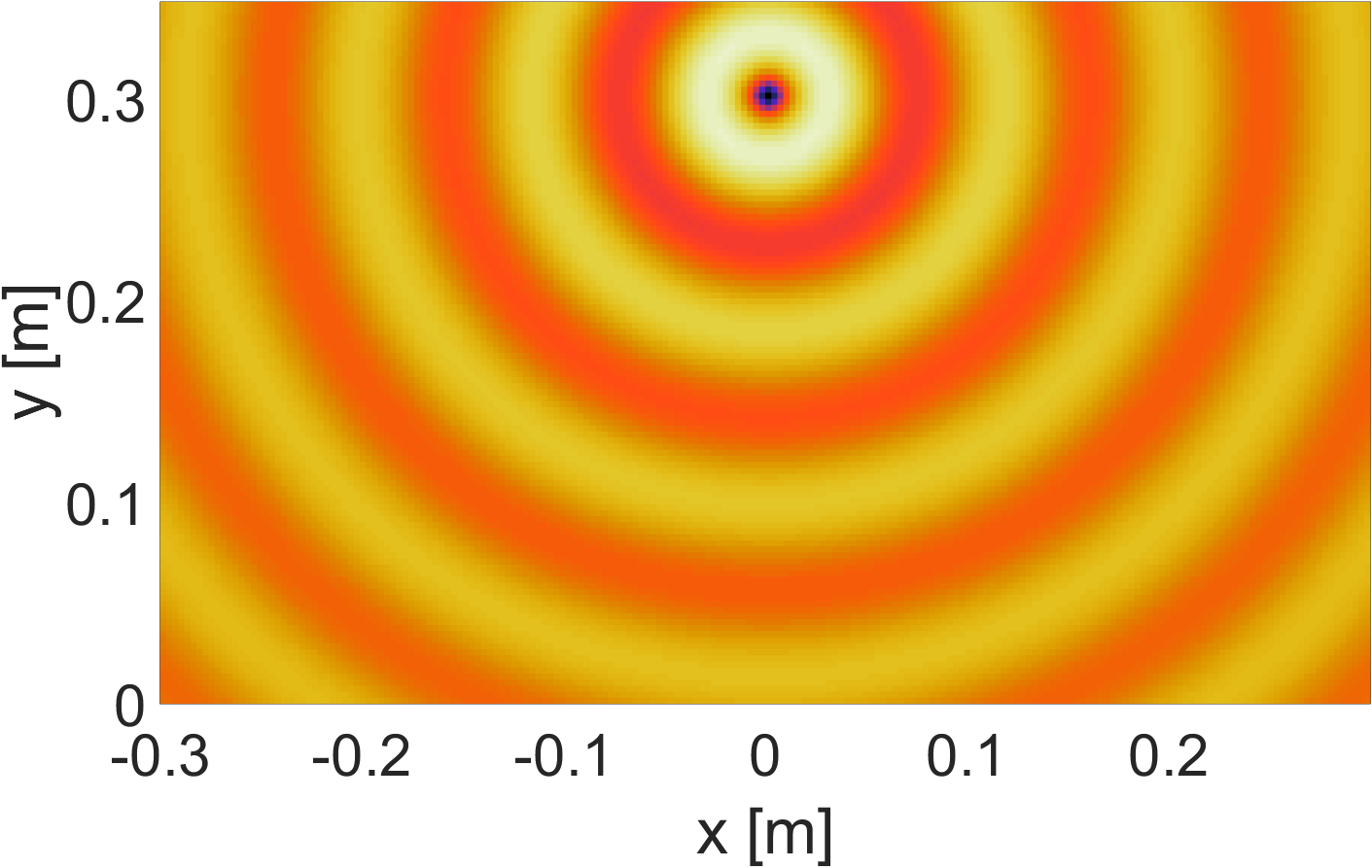}
\caption{}
\label{fig:tx_green_orig}
\end{subfigure}
\begin{subfigure}[h]{\columnwidth}
\centering
\includegraphics[width=5.9cm, clip=true]{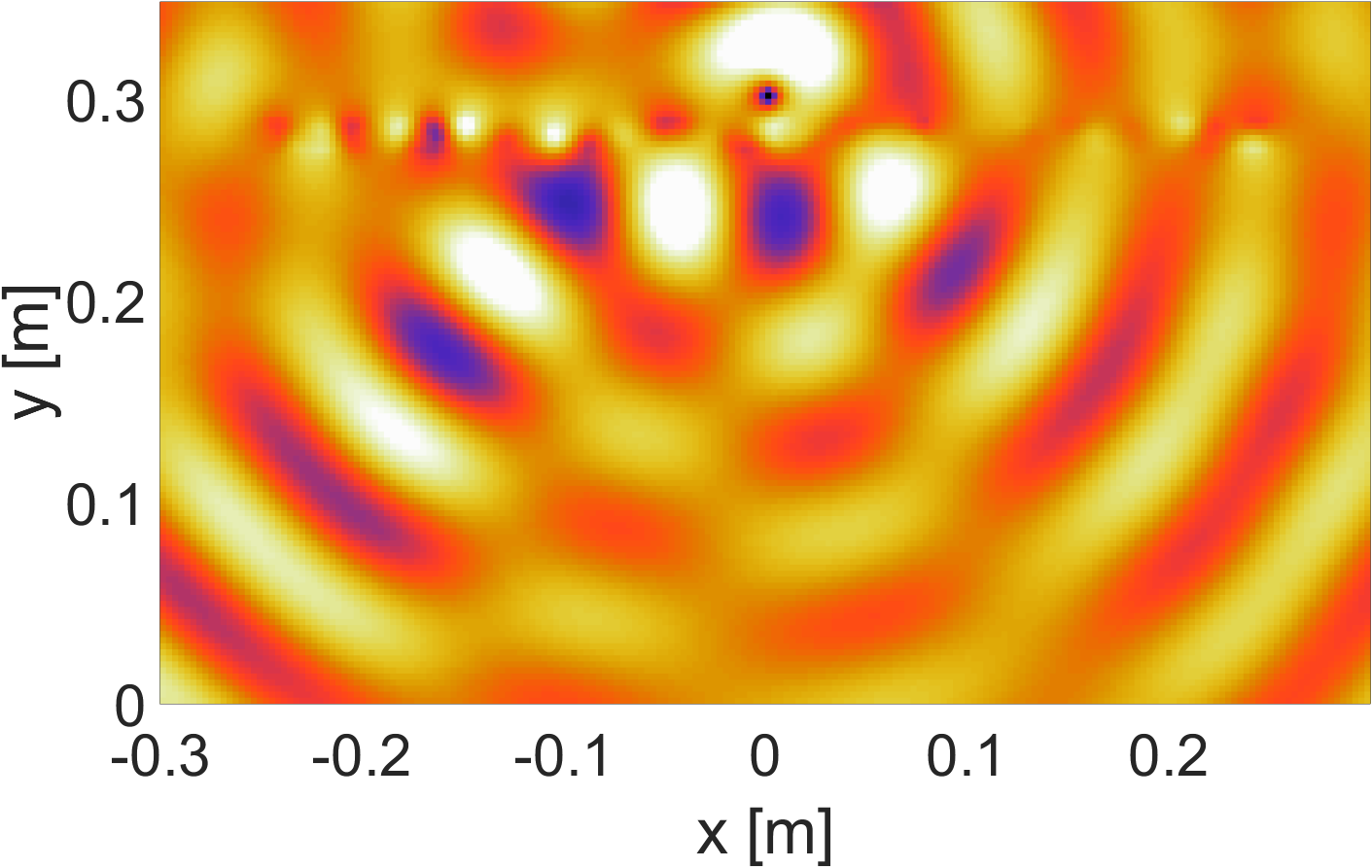}
\caption{}
\label{fig:tx_green_opt}
\end{subfigure}
\caption{\small Green's functions of the compressive antenna operating in transmission mode at $3.5$GHz with (a) the original design and (b) optimized design.}
\end{figure}

Figure \ref{fig:tx_permittivity} displays the optimized permittivity distribution for the transmission mode problem, and Figure \ref{fig:tx_log2} displays the $\log_2$ of the singular values of the Green's function matrices for the original and optimized antennas. In this configuration, the optimized design increases the channel capacity by $40$ bits per pixel and decreases the condition number by a factor of $60$, from approximately $1400$ to $23$. Figures \ref{fig:tx_green_orig} and \ref{fig:tx_green_opt} display the Green's function of the middle antenna radiating at $3.5$GHz for the original and optimized antennas respectively. Without the PEC reflector or any scattering dielectrics, the transmitting antenna radiates isotropically, and so the phase-front is constant as a function of distance from the antenna. In comparison, the phase-front of the optimized design exhibits a pseudo-random pattern, which conveys an increased amount of information compared to the isotropic phase front.


\section{Conclusions}
\label{sec:conclusions}
This paper describes a novel, model-based antenna design method for high-sensing-capacity imaging applications. By optimizing the channel capacity of the dyadic Green's function matrix, the new approach generates antenna configurations with improved sensing and imaging capabilities. Preliminary design results using a 2D FDFD forward model for two antenna configurations, operating in transmission mode and reflection mode, demonstrate how the novel approach generates antenna configurations with improved channel capacity. Although this paper focused on the simplified design approach, the theory of the general design approach allows it to be applied to realistic scenarios. %Furthermore, the flexibility of the generalized design approach allows it to be applied to realistic scenarios.  %Although this paper focused on the simplified design approach, it also prese%Although this paper focused on the simplified design approach, we have presented the theory for the more general design approach and plan to develop it in future work. 

%By optimizing the channel capacity of the dyadic Green's function matrix, the new approach generates antenna configurations with improved sensing and imaging capabilities. Design results for two antenna configurations, which operated in transmission mode and reflection mode, demonstrated how the new method   %We have presented a novel model-based antenna design method for compressive sensing imaging applications. 

\bibliography{SICA-TA}
\bibliographystyle{IEEEtran}

\end{document}